\documentclass[12pt]{article}

\usepackage{amsfonts,amsmath}
\usepackage{latexsym}

\author{ K\'aroly J. B\"or\"oczky\footnote{Supported by
OTKA grants 068398 and 75016, and by the EU Marie Curie TOK
project DiscConvGeo, and FP7 IEF grant GEOSUMSETS.}, Keith M. Ball}
\title{Stability of some versions of the Pr\'ekopa-Leindler inequality}

\newcommand{\proof}{\noindent{\it Proof: }}
\newcommand{\proofbox}{\mbox{ $\Box$}\\}
\newcommand{\R}{\mathbb{R}}

\newtheorem{lemma}{LEMMA}[section]
\newtheorem{theo}[lemma]{THEOREM}

\newtheorem{coro}[lemma]{COROLLARY}

\begin{document}

\maketitle

\begin{abstract}
Two
consequences of the stability version of the
one dimensional Pr\'ekopa-Leindler inequality are presented.
One is the stability version of the Blaschke-Santal\'o
inequality, and the other is
 a stability version of the
Pr\'ekopa-Leindler inequality for even functions
in higher dimensions, where a
recent stability version of the Brunn-Minkowski
inequality is also used in an essential way.
\end{abstract}

\section{The problem}

Our main theme is some consequences of the Pr\'ekopa-Leindler inequality
in one dimension. The inequality itself, due to
A. Pr\'ekopa \cite{Pre71} and L. Leindler \cite{Lei72}, was generalized in
 A. Pr\'ekopa \cite{Pre73} and
\cite{Pre75}, C. Borell \cite{Bor75}, and in H.J. Brascamp, E.H. Lieb \cite{BrL76}.
Various applications are provided and surveyed in K.M. Ball \cite{Bal},
F. Barthe \cite{Bar}, and R.J. Gardner \cite{Gar02}.
The following multiplicative version from \cite{Bal}, is often more useful
and is more convenient for geometric applications.

\begin{theo}[Pr\'ekopa-Leindler]
If $m,f,g$ are non-negative integrable functions on $\R$
satisfying $m(\frac{r+s}2)\geq \sqrt{f(r)g(s)}$ for
$r,s\in\mathbb{R}$, then
$$
\int_{\R}  m\geq \sqrt{\int_{\R}f \cdot \int_{\R}g}.
$$
\end{theo}

S. Dubuc \cite{Dub77} characterized the equality case  if the
integrals of
 $f,g,m$ above are positive, and K.M. Ball, K.J. B\"or\"oczky \cite{BB10}
  even provided the following stability version.

\begin{theo}
\label{PLstab}
There exists an positive absolute constant $c$ with
the following property:
If $m,f,g$ are non-negative integrable functions
with positive integrals on $\R$ such that
$m$ is log-concave,  $m(\frac{r+s}2)\geq \sqrt{f(r)g(s)}$ for
$r,s\in\mathbb{R}$,  and
$$
\int_{\R}  m\leq (1+\varepsilon) \sqrt{\int_{\R}f \cdot \int_{\R}g},
$$
for $\varepsilon>0$, then
there exist $a>0$, $b\in\R$ such that
\begin{eqnarray*}
\int_{\R}|f(t)-a\,m(t+b)|\,dt&\leq &
c\cdot\sqrt[3]{\varepsilon}|\ln \varepsilon|^{\frac43}\cdot
\int_{\R}m(t)\,dt \\
\int_{\R}|g(t)-a^{-1}m(t-b)|\,dt&\leq &
c\cdot\sqrt[3]{\varepsilon}|\ln \varepsilon|^{\frac43}\cdot
\int_{\R}m(t)\,dt.
\end{eqnarray*}
\end{theo}
{\bf Remark } If $f$ and $g$ are log-concave probability distributions
then $a=1$ can be assumed, and if in addition $f$ and $g$ have
the same expectation,
then even $b=0$ can be assumed.\\

As it was observed by C. Borell \cite{Bor75},
and later independently by K.M. Ball \cite{Bal},
assigning to any function $H:[0,\infty]\to[0,\infty]$
the function  $h:\R\to[0,\infty]$ defined by
$h(x)=H(e^x)e^x$,
we have the version Theorem~\ref{PLB} of the Pr\'ekopa-Leindler
inequality. We note that if $H$
is log-concave and decreasing, then $h$
 is log-concave.

\begin{theo}
\label{PLB}
If $M,F,G:[0,\infty]\to[0,\infty]$  integrable functions
satisfy  $M(\sqrt{rs})\geq \sqrt{F(r)G(s)}$ for
$r,s\geq 0$, then
$$
\int_0^\infty  M\geq \sqrt{\int_0^\infty F\cdot \int_0^\infty G}.
$$
\end{theo}

Therefore we deduce the following statement by
Theorem~\ref{PLstab}:

\begin{coro}
\label{PLBstab}
There exists a positive absolute constant $c$ with
the following property:
If $M,F,G:[0,\infty]\to[0,\infty]$ are integrable functions
with positive integrals such that
$M$ is log-concave and decreasing,  $M(\sqrt{rs})\geq \sqrt{F(r)G(s)}$ for
$r,s\in[0,\infty]$,  and
$$
\int_0^\infty  M\leq (1+\varepsilon)
\sqrt{\int_0^\infty F \cdot \int_0^\infty G},
$$
for $\varepsilon>0$, then
there exist $a,b>0$, such that
\begin{eqnarray*}
\int_0^\infty|F(t)-a\,M(b\,t)|\,dt&\leq &
c\cdot\sqrt[3]{\varepsilon}|\ln \varepsilon|^{\frac43}\cdot
\int_0^\infty M(t)\,dt \\
\int_0^\infty|G(t)-a^{-1}M(b^{-1}t)|\,dt&\leq &
c\cdot\sqrt[3]{\varepsilon}|\ln \varepsilon|^{\frac43}\cdot
\int_0^\infty M(t)\,dt.
\end{eqnarray*}
\end{coro}
{\bf Remark }
If in adddition
$F$ and $G$ are decreasing log-concave probability distributions
then $a=b$ can be assumed. The condition that
$M$ is log-concave and decreasing can be replaced by the
one that $M(e^t)$ is log-concave.

\section{A stability version of the Blaschke-Santal\'o
inequality}

Based on the approach in the PhD thesis
K.M. Ball \cite{Bal},
in this section we show how to provide a stability version
of the Blaschke-Santal\'o inequality inequality
using the stability version Corollary~\ref{PLBstab} of the
Pr\'ekopa-Leindler inequality in one one dimension.

We write $o$ to denote the origin of $\R^n$,
$\langle \cdot,\cdot\rangle$ to denote
the standard scalar product.
We write $|\cdot|$ to denote
the Lebesgue measure
in $\R^n$, where the Lebesgue measure of the empty set is $0$.
Let $B^n$ be the unit Euclidean ball
with volume $\kappa_n=|B^n|$.
A convex body $K$ in $\R^n$ is a compact convex set
with non--empty interior.
If $z\in{\rm int}K$,
then the polar of $K$ with respect to $z$ is the convex body
$$
K^z=\{x\in\R^n:\,\langle x-z,y-z\rangle\leq 1\mbox{ for any $y\in K$}\}.
$$
It is easy to see that $(K^z)^z=K$, and the volume product
$|K|\cdot|K^z|$ is affine invariant.
According to L.A. Santal\'o \cite{San49}
 (see also M. Meyer and A. Pajor \cite{MeP90}), there exists a unique
$z\in{\rm int}K$ minimizing
$|K^z|$, which is called the Santal\'o point
of $K$. In this case $z$ is the centroid of $K^z$.
The celebrated  Blaschke-Santal\'o
inequality states that if $z$ is the  Santal\'o point
(or centroid) of $K$, then
\begin{equation}
\label{BS}
|K|\cdot |K^z|\leq\kappa_n^2,
\end{equation}
with equality if and only if $K$ is an ellipsoid.
The inequality was proved by W. Blaschke \cite{Bla17} for $n\leq 3$,
 and by L.A. Santal\'o \cite{San49} for all $n$.
The case of equality was characterized by
J. Saint-Raymond \cite{Sai81} among $o$-symmetric convex bodies,
and by C.M. Petty \cite{Pet85} among all convex bodies
(see also M. Meyer and A. Pajor \cite{MeP90},
D. Hug \cite{Hug96}, and M. Meyer and S. Reisner \cite{MeR06}
 for simpler proofs).

A natural tool is the  Banach-Mazur distance $\delta_{\rm BM}(K,M)$
of the convex bodies $K$ and $M$, which is defined by
$$
\delta_{\rm BM}(K,M)=\ln\min\{\lambda\geq 1:\,
K-x\subset \Phi(M-y)\subset \lambda(K-x)\mbox{ for }
\Phi\in{\rm GL}(n),x,y\in\R^n\}.
$$
In particular, if $K$ and $M$ are $o$-symmetric, then
$x=y=o$ can be assumed, and in this case
$\delta_{\rm BM}(K,M)=\delta_{\rm BM}(K^o,M^o)$. It follows
from a theorem of F. John \cite{Joh37}
that  $\delta_{\rm BM}(K,B^n)\leq \ln n$
for any convex body $K$ in $\R^n$
(see also K.M. Ball \cite{Bal97}).

K.J. B\"or\"oczky \cite{Bor} proved a stability
version of the Blaschke-Santal\'o
inequality. One of the main tool in that paper is to reduce the
problem to $o$-symmetric convex
bodies with axial rotational symmetry; namely, combining
Theorem~1.4 and Lemma~2.1 in \cite{Bor} yields the following.

\begin{lemma}
\label{rounding}
For any $n\geq 2$ there exists $\tilde{\gamma}>0$ depending, such that
if $K$ is a convex body  in $\R^n$  with Santal\'o point
$z$, then one finds an $o$-symmetric convex body $C$ with
axial rotational symmetry, and
satisfying $\delta_{\rm BM}(C,B^n)\geq \tilde{\gamma}\delta_{\rm BM}(K,B^n)^2$
and $|C|\cdot |C^o|\geq |K|\cdot |K^z|$.
\end{lemma}
{\bf Remark: } If $K$ is $o$-symmetric, then even
$\delta_{\rm BM}(C,B^n)\geq \tilde{\gamma}\delta_{\rm BM}(K,B^n)$.\\

Now we are ready to prove our main result in this section:

\begin{theo}
\label{BSstab}
If $K$ is a  convex body in $\R^n$, $n\geq 3$, with Santal\'o point $z$,
and
$$
(1+\varepsilon)|K|\cdot|K^z|>\kappa_n^2
\mbox{ \ for $\varepsilon>0$,}
$$
then for some $\gamma>0$ depending only on $n$, we have
$$
\delta_{\rm BM}(K,B^n)<\gamma\,\varepsilon^{\frac1{3(n+1)}}
|\log \varepsilon|^{\frac4{3(n+1)}}.
$$
\end{theo}
{\bf Remark: } If $K$ is
$o$-symmetric, then  the exponent $\frac1{3(n+1)}$ occuring
in Theorem~\ref{BSstab} can be replaced by $\frac2{3(n+1)}$.

Taking $K$ to be the convex body resulting from $B^n$
by cutting off two opposite caps of volume $\varepsilon$  shows that
 the exponent $1/(3(n+1))$ cannot be replaced by anything larger
than $2/(n+1)$ even for $o$-symmetric
convex bodies with axial rotational
symmetry. Therefore the exponent of $\varepsilon$
is of the correct order. In addition if the error in
Corollary~\ref{PLBstab} is reduced to $\varepsilon$, then
we would have the the stability
version of the Blaschke-Santal\'o
inequality of the correct order.

We note that
the exponent of $\varepsilon$ is $1/(6n)$
in the stability version of the Blaschke-Santal\'o
inequality proved in K.J. B\"or\"oczky~\cite{Bor}.\\

\noindent{\it Proof of Theorem~\ref{BSstab}: }
Let $C$ be the $o$-symmetric convex
body provided by Lemma~\ref{rounding}, let $u$ be
a unit vector and $\alpha>0$ such that $\alpha u\in\partial C$,
and a section
$C\cap(u^\bot+tu)$ for $t\in (-\alpha,\alpha)$ is
an $(n-1)$-ball of radius $\varphi(t)$ and
of area $F(t)=\varphi(t)^{n-1}\kappa_{n-1}$.
In turn $\alpha^{-1} u\in\partial C^o$,
and if $t\in (-\alpha^{-1},\alpha^{-1})$, then
$C^o\cap(u^\bot+tu)$ is
an $(n-1)$-ball of radius $\psi(t)$ and
of area $G(t)=\psi(t)^{n-1}\kappa_{n-1}$. We observe that for $t\in (-1,1)$,
$B^n\cap(u^\bot+tu)$ is
an $(n-1)$-ball of radius $(1-t^2)^{\frac12}$ and
of area $M(t)=(1-t^2)^{\frac{n-1}2}\kappa_{n-1}$.
We define $F(t)=0$, $G(t)=0$ and $M(t)=0$
if $t\geq\alpha$, $t\geq\alpha^{-1}$,
and $t\geq 1$, respectively.
For $v\in u^\bot$, $r\in (-\alpha,\alpha)$
and $s\in (-\alpha^{-1},\alpha^{-1})$, we have
$$
\varphi(r)\cdot\psi(s)=\langle ru+\varphi(r)v,su+\psi(s)v\rangle
-rs\leq 1-rs.
$$
In particular $M(\sqrt{rs})\geq \sqrt{F(r)G(s)}$, and
$$
\left(\int_0^\infty M\right)^2=\frac{\kappa_n^2}4\leq
\frac{(1+\varepsilon)^{-1}|C|\cdot|C^0|}4=
(1+\varepsilon)\left(\int_{0}^\infty F\right)
\left(\int_{0}^\infty G\right).
$$
Therefore we may apply Corollary~\ref{PLBstab}, and
deduce that there exist $a,b>0$ such that
\begin{equation}
\label{BSFM}
\int_0^\infty|a\,F(b\,t)-M(t)|\,dt\leq
c\cdot\sqrt[3]{\varepsilon}|\ln \varepsilon|^{\frac43}\cdot
\kappa_n.
\end{equation}
Let $\Phi$ be the linear transform such that
$\Phi u=b^{-1}u$, and if $v\in u^\bot$ then
$\Phi v=a^{\frac1{n-1}}v$. Therefore
$\widetilde{C}=\Phi C$ is an $o$-symmetric
convex body with axial symmetry around $\R u$
such that
\begin{equation}
\label{BMCK}
\delta_{\rm BM}(\widetilde{C},B^n)=
\delta_{\rm BM}(C,B^n)\geq \tilde{\gamma}\delta_{\rm BM}(K,B^n)^2,
\end{equation}
and  the area of
$\widetilde{C}\cap(u^\bot+tu)$ is $a\,F(b\,t)$
for any $t\in (0,b^{-1}\alpha)$.
In particular
$$
|\widetilde{C}\Delta B^n|\leq
2c\cdot\sqrt[3]{\varepsilon}|\ln \varepsilon|^{\frac43}\cdot
\kappa_n.
$$
Since there exists some $\gamma_0>0$ depending on $n$ such that
$$
(1-\gamma_0|\widetilde{C}\Delta B^n|^{\frac2{n+1}})B^n
\subset \widetilde{C}\subset
(1+\gamma_0|\widetilde{C}\Delta B^n|^{\frac2{n+1}})B^n,
$$
we conclude Theorem~\ref{BSstab} by (\ref{BMCK}).
\proofbox

As it is explained in
K.J. B\"or\"oczky~\cite{Bor}, the stability
version Theorem~\ref{BSstab} yields
stability versions with the same order of the error term
for two basic affine invariant inequalities.
The first is the affine isoperimetric inequality
of W. Blaschke \cite{Bla16} or \cite{Bla22}
(see L.A. Santal\'o \cite{San49} for $n\geq 4$),
and the other is the isoperimetric inequality
for the geominimal surface area by
C.M. Petty \cite{Pet85}. The
monograph K. Leichtwei{\ss} \cite{Lei98}, and
and the survey paper E. Lutwak \cite{Lut93}
provide introduction into
the by now classical theory of these notions.

\section{The stability version of the Brunn-Minkowski
inequality due to Figalli, Maggi, Pratelli}

 For any $\alpha,\beta>0$, and
measurable sets $X,Y,Z\subset\R^n$ with
$$
\alpha X+\beta Y\subset Z,
$$
 the Brunn-Minkowski inequality says that
$$
|Z|^{\frac1n}\geq \alpha |X|^{\frac1n}+\beta |Y|^{\frac1n}.
$$
Obviously the case $\alpha=\beta=1$ yields the general case.

Since the days of H. Minkowski, there are stability versions
of the Brunn-Minkowski inequality if $X$ and $Y$ convex bodies, mostly
in terms of the so called Hausdorff metric,
see the survey paper H. Groemer \cite{Gro93}. In higher
dimensions, the best estimates are due to
V.I. Diskant \cite{Dis73} and H. Groemer \cite {Gro88}.

Recently A. Figalli, F. Maggi, A. Pratelli \cite{FMP1} and \cite{FMP1}
obtained an optimal stability version of the Brunn-Minkowski
inequality in terms of the volume difference.
To define the ``homothetic distance'' $A(K,C)$
of convex bodies $K$ and $C$, let $\alpha=|K|^{\frac{-1}n}$ and
$\beta=|C|^{\frac{-1}n}$, and let
$$
A(K,C)=\min\left\{|\alpha K\Delta (x+\beta C)|:\,x\in\R^n\right\}.
$$
We observe that $|\alpha K\cap(x+\beta C)|^{\frac1n}$ is a concave
function of $x\in \alpha K -\beta C$
by the Brunn-Minkowski inequality. Therefore if both
$K$ and $C$ are $o$-symmetric, and $|C|=|K|$, then
\begin{equation}
\label{Aosymm}
A(K,C)=|K\Delta C|/|K|.
\end{equation}
Next let
$$
\sigma(K,C)=\max\left\{\frac{|C|}{|K|},\frac{|K|}{|C|}\right\}.
$$

\begin{theo}[Figalli,Maggi,Pratelli]
\label{Maggi}
For $\gamma^*=(\frac{(2-2^{\frac{n-1}{n}})^{\frac32}}{122n^7})^2$,
and  any convex bodies $K$ and $C$ in $\R^n$,
$$
|K+C|^{\frac1n}\geq (|K|^{\frac1n}+|C|^{\frac1n})
\left[1+\frac{\gamma^*}{\sigma(K,C)^{\frac1n}}\cdot A(K,C)^2\right].
$$
\end{theo}

We will need the product form of the Brunn-Minkowski inequality.
Since
\begin{eqnarray*}
\frac12\left(|K|^{\frac1n}+|C|^{\frac1n}\right)&=&
|K|^{\frac1{2n}}|C|^{\frac1{2n}}
\left[1+\frac12
\left(\sigma(K,C)^{\frac1{4n}}-\sigma(K,C)^{\frac{-1}{4n}}\right)^2\right] \\
&\geq &|K|^{\frac1{2n}}|C|^{\frac1{2n}}
\left[1+\frac{(\sigma(K,C)-1)^2}{32n^2\sigma(K,C)^{\frac{4n-1}{2n}}}\right],
\end{eqnarray*}
we conclude with $\sigma=\sigma(K,C)$ that
\begin{equation}
\label{Maggiprod}
\left|\mbox{$\frac12$}(K+C)\right|\geq \sqrt{|K|\cdot|C|}\left[1+
\frac{(\sigma-1)^2}{32n\sigma^2}
+\frac{n\gamma^*}{\sigma^{\frac1n}}\cdot A(K,C)^2\right].
\end{equation}

\section{Pr\'ekopa-Leindler inequality
in higher dimensions for even functions}

Let $f,g,m:\R^n\to[0,\infty]$ such that
$m(\frac{x+y}2)\geq \sqrt{f(x)g(y)}$ for
$r,s\in\mathbb{R}^n$, and for $t>0$, let
\begin{eqnarray*}
\Phi_t&=&
\{x\in\R^n:\,f(x)\geq t\} \mbox{ \ and \ } F(t)=|\Phi_t|\\
\Psi_t&=&\{x\in\R^n:\,g(x)\geq t\} \mbox{ \ and \ } G(t)=|\Psi_t|\\
\Omega_t&=&\{x\in\R^n:\,m(x)\geq t\}\mbox{ \ and \ } M(t)=|\Omega_t|.
\end{eqnarray*}
As it was observed in K.M. Ball \cite{Bal}, the condition
on $f,g,m$ yields that if  $\Phi_r,\Psi_s\neq\emptyset$ for $r,s>0$, then
\begin{equation}
\label{sections}
\mbox{$\frac12$}(\Phi_r+\Psi_s)\subset \Omega_{\sqrt{rs}}.
\end{equation}
Therefore the Brunn-Minkowski inequality yields that
\begin{equation}
\label{minksum}
M(\sqrt{rs})\geq \left(\frac{F(r)^{\frac1n}+G(s)^{\frac1n}}2\right)^n
\geq\sqrt{F(r)\cdot G(s)}
\end{equation}
for all $r,s>0$.
In particular we deduce the Pr\'ekopa-Leindler inequality
by Theorem~\ref{PLB}, as
$$
\int_{\R^n}m=\int_0^\infty M(t)\,dt\geq
\sqrt{\int_0^\infty F(t)\,dt \cdot \int_0^\infty G(t)\,dt}
=\sqrt{\int_{\R^n}f \cdot \int_{\R^n}g}.
$$

The main goal of this section is to prove
a stability version of Pr\'ekopa-Leindler inequality
for at least for even functions. First let
$$
\omega(\varepsilon)=\sqrt[3]{\varepsilon}|\ln \varepsilon|^{\frac43},
$$
which is the error estimate in Theorem~\ref{PLstab}
(and hence the error estimate in Corollary~\ref{PLBstab}).
If $\varphi$ and $\psi$ are real functions, then
we write $\varphi\ll\psi$ if there exists a $\gamma>0$
depending only on $n$ such that $|\varphi|\leq \gamma\cdot \psi$.

\begin{theo}
\label{PLhstab}
If $m,f,g:\,\R^n\to [0,\infty]$ are even and integrable  such that
$m$ is log-concave,  $m(\frac{x+y}2)\geq \sqrt{f(x)g(y)}$ for
$x,y\in\mathbb{R}^n$,  and
$$
\int_{\R^n}  m\leq (1+\varepsilon) \sqrt{\int_{\R^n}f \cdot \int_{\R^n}g}
$$
for $\varepsilon>0$, then
there exist $a>0$ such that
\begin{eqnarray*}
\int_{\R^n}|a\,f(x)-m(x)|\,dx&\ll &\sqrt{\omega(\varepsilon)}
\int_{\R^n}m \\
\int_{\R^n}|a^{-1}g(x)-m(x)|\,dx&\ll &\sqrt{\omega(\varepsilon)}
\int_{\R^n}m.
\end{eqnarray*}
\end{theo}
\proof As in the one dimensional case, we may assume
that $f,g:\,\R^n\to [0,\infty]$ are even and log-concave probability
distributions.
 We may also assume that
$\varepsilon\in(0,\varepsilon_0)$ where $\varepsilon_0\in(0,1)$
is chosen in a  suitable way and depends only on $n$.

We define $\Phi_t,\Psi_t,\Omega_t$ and $F(t),G(t),M(t)$
analogously as at the beginning of the section.
We observe that $\Phi_t,\Psi_t,\Omega_t$
are $o$-symmetric convex bodies, and
$F(t),G(t),M(t)$ are decreasing and log-concave,
and $F,G$ are probability distributions
on $[0,\infty]$. Since $\int_0^\infty M=\int_{\R^n}  m\leq (1+\varepsilon)$,
it follows from Corollary~\ref{PLBstab} that
there exists some $b>0$ such that
\begin{eqnarray*}
\int_0^\infty|bF(bt)-M(t)|\,dt&\ll &\omega(\varepsilon)\\
\int_0^\infty|b^{-1}G(b^{-1}t)-M(t)|\,dt&\ll &\omega(\varepsilon).
\end{eqnarray*}
We may assume that $b\geq 1$. For $x\in \R^n$, we define
\begin{eqnarray*}
\tilde{f}(x)&= &b^{-1}f(b^{\frac{-1}n}x)\\
\tilde{g}(x)&= &b\,g(b^{\frac{1}n}x).
\end{eqnarray*}

The main strategy of the proof is as follows.
First we verify
\begin{equation}
\label{tildefgdist}
\int_{\R^n}|\tilde{f}(x)-\tilde{g}(x)|\,dx\ll \sqrt{\omega(\varepsilon)}.
\end{equation}
Along the way, we establish $b-1\ll\sqrt{\omega(\varepsilon)}$,
which in turn yields
\begin{equation}
\label{fgdist}
\int_{\R^n}|f(x)-g(x)|\,dx\ll \sqrt{\omega(\varepsilon)}.
\end{equation}
Finally we conclude Theorem~\ref{PLhstab} from
(\ref{fgdist}) and (\ref{sections}).

For $t>0$, let
\begin{eqnarray*}
\widetilde{\Phi}_t&=&
\{x\in\R^n:\,\tilde{f}(x)\geq t\}
\mbox{ \ \ \ where $\widetilde{\Phi}_t=b^{\frac{1}n}\Phi_{bt}$
if $\widetilde{\Phi}_t\neq\emptyset$} \\
\widetilde{\Psi}_t&=&\{x\in\R^n:\,\tilde{g}(x)\geq t\}
\mbox{ \ \ \ where $\widetilde{\Psi}_t=b^{\frac{-1}n}\Psi_{b^{-1}t}$
if $\widetilde{\Psi}_t\neq\emptyset$}.
\end{eqnarray*}
These sets satisfy
\begin{eqnarray*}
\int_0^\infty|\,|\widetilde{\Phi}_t|-M(t)|\,dt&\ll &\omega(\varepsilon)\\
\int_0^\infty|\,|\widetilde{\Psi}_t|-M(t)|\,dt&\ll &
\omega(\varepsilon),
\end{eqnarray*}
and (\ref{sections}) yields that
if $\widetilde{\Phi}_t\neq\emptyset$ and $\widetilde{\Psi}_t\neq\emptyset$
for $t>0$, then
\begin{equation}
\label{widetildeOmega}
 \mbox{$\frac12$}(b^{\frac{-1}n}\widetilde{\Phi}_t
+b^{\frac{1}n}\widetilde{\Psi}_t)\subset\Omega_t.
\end{equation}
The main task is estimate the $L_1$ distance
of $\tilde{f}$ and $\tilde{g}$ using
$$
\int_{\R^n}|\tilde{f}(x)-\tilde{g}(x)|\,dx=
\int_0^\infty|\widetilde{\Phi}_t\Delta\widetilde{\Psi}_t|\,dt.
$$
We dissect $[0,\infty)$ into $I$ and $J$,
where $t\in I$, if $\frac34\,M(t)<|\widetilde{\Phi}_t|<\frac54\,M(t)$
and $\frac34\,M(t)<|\widetilde{\Psi}_t|<\frac54\,M(t)$, and
$t\in J$ otherwise. If $t\in J$, then
$$
|\widetilde{\Phi}_t\Delta\widetilde{\Psi}_t|\leq |\widetilde{\Phi}_t|+
|\widetilde{\Psi}_t|\leq 10\left(|\,|\widetilde{\Phi}_t|-M(t)|+
|\,|\widetilde{\Psi}_t|-M(t)|\right).
$$
Therefore
\begin{equation}
\label{tildefgJ}
\int_J|\widetilde{\Phi}_t\Delta\widetilde{\Psi}_t|\,dt\ll
\omega(\varepsilon).
\end{equation}
In addition if $\varepsilon_0$ is small enough, then
\begin{equation}
\label{Jsize}
\int_J M(t)\,dt\leq 4\int_J\left(|\,|\widetilde{\Phi}_t|-M(t)|+
|\,|\widetilde{\Psi}_t|-M(t)|\right)\,dt
\ll\omega(\varepsilon)<\frac12.
\end{equation}

Turning to $I$, it follows form the Pr\'ekopa-Leindler inequality
and (\ref{Jsize}) that
\begin{equation}
\label{Isize}
\int_I M(t)\,dt\geq 1-\int_J M(t)\,dt>\frac12.
\end{equation}
For $t\in I$, we define
 $\alpha(t)=|\widetilde{\Phi}_t|/M(t)$ and
$\beta(t)=|\widetilde{\Psi}_t|/M(t)$, and hence
$\frac34<\alpha(t),\beta(t)<\frac54$, and
\begin{equation}
\label{alphabeta}
\int_0^\infty M(t)\cdot\left(|\alpha(t)-1|+|\beta(t)-1|\right)\,dt\ll
\omega(\varepsilon).
\end{equation}
In addition let
\begin{eqnarray*}
\sigma(t)&=&\sigma\left(b^{\frac{-1}n}\widetilde{\Phi}_t,
b^{\frac{1}n}\widetilde{\Psi}_t\right)=
\max\left\{\frac{b^2\beta(t)}{\alpha(t)},
\frac{\alpha(t)}{b^2\beta(t)}\right\}\\
\eta(t)&=&
\frac{(\sigma(t)-1)^2}{32n\sigma(t)^2}
+\frac{n\gamma^*}{\sigma(t)^{\frac1n}}\cdot
A(\widetilde{\Phi}_t,\widetilde{\Psi}_t)^2,
\end{eqnarray*}
where $\gamma^*$ comes from Theorem~\ref{Maggi}.
It follows from $\alpha(t),\beta(t)>\frac34$,
(\ref{Maggiprod}) and (\ref{widetildeOmega})
that
\begin{eqnarray*}
M(t)&\geq& M(t)\cdot\sqrt{\alpha(t)\cdot \beta(t)}(1+\eta(t))\\
&\geq& M(t)\cdot\left(1-\max\{0,1-\alpha(t)\}-\max\{0,1-\beta(t)\}\right)
(1+\eta(t))\\
&\geq& M(t)\cdot(1-|\alpha(t)-1|-|\beta(t)-1|+\mbox{$\frac12$}\,\eta(t)).
\end{eqnarray*}
In particular (\ref{alphabeta}) yields
\begin{equation}
\label{Ieta}
\int_I  M(t)\cdot\eta(t)\,dt\ll \omega(\varepsilon).
\end{equation}

Next we estimate $b$. Let $t\in I$.
If $\alpha(t)\geq b\beta(t)$ then
$$
|\alpha(t)-1|+|\beta(t)-1|\geq\frac{\sqrt{b}-1}{\sqrt{b}}\geq
\frac{b-1}{2b}\geq\frac{(b-1)^2}{32n b^2}.
$$
If $\alpha(t)< b\beta(t)$ then $\sigma(t)> b$, and
$$
\eta(t)>\frac{(b-1)^2}{32n b^2}.
$$
We deduce by (\ref{Isize}), (\ref{alphabeta})
 and (\ref{Ieta}) that
$$
\frac{(b-1)^2}{64n b^2}\leq \int_I M(t)\cdot\frac{(b-1)^2}{32n b^2}\,dt
\leq \int_0^\infty M(t)\cdot
\left(\eta(t)+|\alpha(t)-1|+|\beta(t)-1|\right)\,dt\ll
\omega(\varepsilon).
$$
Since $\frac{b-1}{b}>\frac12$ if $b>2$, we deduce that
\begin{equation}
\label{best}
b-1\ll\sqrt{\omega(\varepsilon)}.
\end{equation}
It also follows that $\sigma(t)<2$ if $\varepsilon_0$ is small
enough, and hence $A(\widetilde{\Phi}_t,\widetilde{\Psi}_t)^2\ll\eta(t)$.

For $t\in I$, we deduce using (\ref{Aosymm}) that
\begin{eqnarray*}
|\widetilde{\Phi}_t\Delta\widetilde{\Psi}_t|^2&\leq&
3\left[|\alpha(t)^{\frac{-1}n}\widetilde{\Phi}_t\Delta
\beta(t)^{\frac{-1}n}\widetilde{\Psi}_t|^2+
|\alpha(t)^{\frac{-1}n}\widetilde{\Phi}_t\Delta\widetilde{\Phi}_t|^2
+|\beta(t)^{\frac{-1}n}\widetilde{\Psi}_t\Delta\widetilde{\Psi}_t|^2
\right]\\
&=&3\left[A(\widetilde{\Phi}_t,\widetilde{\Psi}_t)^2
+|\alpha(t)-1|^2+|\beta(t)-1|^2\right]\cdot M(t)^2.
\end{eqnarray*}
In turn we have
\begin{eqnarray*}
\left(\int_I|\widetilde{\Phi}_t\Delta\widetilde{\Psi}_t|\,dt\right)^2&\leq&
\int_I \frac{|\widetilde{\Phi}_t\Delta\widetilde{\Psi}_t|^2}{M(t)}\,dt
\int_I M(t)\,dt  \\
&\ll& \int_I A(\widetilde{\Phi}_t,\widetilde{\Psi}_t)^2\cdot M(t)\,dt+
\omega(\varepsilon)\\
&\ll&
\int_I \eta(t)\cdot M(t)\,dt+
\omega(\varepsilon)\ll \omega(\varepsilon).
\end{eqnarray*}
Combining this estimate with (\ref{tildefgJ}) yields
\begin{equation}
\label{tildefg}
\int_{\R^n}|\tilde{f}(x)-\tilde{g}(x)|\,dx=
\int_0^\infty|\widetilde{\Phi}_t\Delta\widetilde{\Psi}_t|\,dt
\ll \sqrt{\omega(\varepsilon)}.
\end{equation}

Turning to the $L_1$ distance of $f$ and $g$,
$f(b^{\frac{-1}n}x)\geq f(x)$ for $x\in\R^n$ and
the estimate (\ref{best}) on $b$ yield
\begin{eqnarray*}
\int_{\R^n}|f-\tilde{f}|&\leq&
\int_{\R^n}(f-b^{-1}f)+
\int_{\R^n}[b^{-1}f(b^{\frac{-1}n}x)-b^{-1}f(x)]\,dx\\
&=&2(1-b^{-1})\ll \sqrt{\omega(\varepsilon)}.
\end{eqnarray*}
Similarly $\int_{\R^n}|g(x)-\tilde{g}(x)|\,dx\ll \sqrt{\omega(\varepsilon)}$,
therefore (\ref{tildefg}) implies
\begin{equation}
\label{fg}
\int_{\R^n}|f(x)-g(x)|\,dx=
\int_0^\infty|\Phi_t\Delta\Psi_t|\,dt
\leq \gamma_1\sqrt{\omega(\varepsilon)}
\end{equation}
for a $\gamma_1>0$ depending only on $n$.

Finally, we compare $f$ and $m$.
It follows by (\ref{sections}) that
if  $\Phi_t,\Psi_t\neq\emptyset$ for $t>0$, then
\begin{equation}
\label{sectioncap}
\Phi_t\cap\Psi_t\subset \Omega_t.
\end{equation}
Using the $\gamma_1$ of (\ref{fg}), we have
\begin{eqnarray*}
1-\gamma_1\sqrt{\omega(\varepsilon)}&\leq&
\int_0^\infty|\Phi_t|\,dt-
\int_0^\infty|\Phi_t\Delta\Psi_t|\,dt\\
&\leq&
\int_0^\infty|\Phi_t\cap\Psi_t|\,dt\leq
\int_0^\infty|\Omega_t|\,dt\leq 1+\varepsilon.
\end{eqnarray*}
Since (\ref{sectioncap}) yields that
$$
|\Phi_t\Delta\Omega_t|\leq
|\Phi_t\Delta\Psi_t|+|\Omega_t|-|\Phi_t\cap\Psi_t|,
$$
we conclude that
$$
\int_{\R^n}|f(x)-m(x)|\,dx=
\int_0^\infty|\Phi_t\Delta\Omega_t|\,dt\ll \sqrt{\omega(\varepsilon)}.
$$
Similarly we have $\int_{\R^n}|g(x)-m(x)|\,dx\ll \sqrt{\omega(\varepsilon)}$.
\proofbox

\end{document}